\documentclass[reqno]{amsart}
\usepackage{amsmath}
\usepackage{xspace}
\usepackage{alltt}
\usepackage{tables}
\theoremstyle{definition}

\theoremstyle{plain}
\newtheorem{thm}{Theorem}

\numberwithin{equation}{section}
\allowdisplaybreaks[1]

\newcommand{\qBin}[2]{\genfrac{[}{]}{0pt}{0}{#1}{#2}}

\newcommand{\qhg}{$q$\nobreakdash-hyper\-geometric\xspace}

\newlength{\axellength}

\def\eqn#1{\eqref{eq:#1}}
\def\mycap#1#2{\hbox to \textwidth{\hfil{\textsc{Table #1}.{\ #2}}\hfil}}
\def\tabref#1{\textsc{Table #1}}

\begin{document}


 \makeatletter
 \def\print@eqnum{\tagform@{\normalsize\theequation}}
 \def\endmathdisplay@a{%
   \if@eqnsw \gdef\df@tag{\tagform@{\normalsize\theequation}}\fi
   \if@fleqn \@xp\endmathdisplay@fleqn
   \else \ifx\df@tag\@empty \else \veqno \alt@tag \df@tag \fi
     \ifx\df@label\@empty \else \@xp\ltx@label\@xp{\df@label}\fi
   \fi
   \ifnum\dspbrk@lvl>\m@ne
     \postdisplaypenalty -\@getpen\dspbrk@lvl
     \global\dspbrk@lvl\m@ne
   \fi
 }
 \makeatother

\title[\tiny{A new four parameter $q-$series identity and its partition implications}]
{A new four parameter $q-$series identity  and its partition implications}
\author{Krishnaswami Alladi,  George E. Andrews, and 
Alexander Berkovich}

\address{Department of Mathematics, The University of Florida,
         Gainesville, FL~32611, USA}
\email{alladi@math.ufl.edu}

\address{Department of Mathematics, The Pennsylvania State University,
         University Park, PA~16802, USA}
\email{andrews@math.psu.edu}

\address{Department of Mathematics, The University of Florida,
         Gainesville, FL~32611, USA}
\email{alexb@math.ufl.edu}

\thanks{Research of the first two authors is supported in part by grants from the National 
Science Foundation, and of the third author by a grant from the Number Theory Foundation}

\maketitle

\begin{abstract}
We prove a new four parameter \qhg series identity from which 
the three parameter identity for the G\"ollnitz theorem due to Alladi, 
Andrews, and Gordon follows as a special case by setting one of the 
parameters equal to $0$.  The new identity is equivalent to a  four 
parameter partition theorem which extends the deep theorem of G\"ollnitz 
and thereby settles a problem raised by Andrews thirty years ago. Some 
consequences including a quadruple product extension of Jacobi's triple 
product identity, and prospects of future research are briefly discussed.
\end{abstract}


\section{\bf Introduction}

\medskip

One of the most fundamental results in the theory of partitions is the 1926 
theorem of Schur \cite{26} which we state in the following form:
\begin{thm}(Schur)\\

Let $P_1(n)$ denote the number of partitions of $n$ into distinct parts 
$\equiv-2^1,-2^0\pmod{3}$.

Let $G_1(n)$ denote the number of partitions of $n$ such that the difference between
the parts $\ge3$, with equality only if a part is $\equiv-2^1,-2^0\pmod{3}$. Then
$$
G_1(n)=P_1(n).
$$
\label{thm:1}
\end{thm}
In 1967, G\"ollnitz \cite{25} established a deep partition theorem which we 
prefer to state as follows:
\begin{thm}(G\"ollnitz)\\

Let $P_2(n)$ denote the number of partitions of $n$ into distinct parts 
$\equiv-2^2,-2^1,-2^0\pmod{6}$.

Let $G_2(n)$ denote the  number of partitions of $n$ into parts $\neq  1$ or $3$, 
such that the difference between the parts $\ge6$, with equality only if a part is $\equiv-2^2,-2^1,
-2^0\pmod{6}$. Then
$$
G_2(n)=P_2(n).
$$
\label{thm:2}
\end{thm}
Although Theorem~\ref{thm:1} is not a special case of Theorem~\ref{thm:2}, the result of 
G\"ollnitz can be viewed as the next higher level extension of Schur's 
theorem especially in the forms in which we have stated these 
results.  While attempting to find new partition identities via a computer 
search in 1971, Andrews \cite{19} raised the question whether there exists a 
partition theorem that goes beyond Theorem~\ref{thm:2} in the sense that Theorem~\ref{thm:2} 
may be viewed as going beyond Theorem~\ref{thm:1}.  We settle this problem in the 
affirmative by proving the following result along with its refinement and 
generalization (Theorem~\ref{thm:6}) stated in \S3:
\begin{thm} 

Let $P_3(n)$ denote the number of partitions of $n$ into distinct parts 
$\equiv-2^3,-2^2,-2^1,
\\-2^0\pmod{15}$.

Let $G_3(n)$ denote the number of partitions of $n$ into parts 
$\not\equiv 2^3,2^2,2^1,2^0\pmod{15}$, such that the difference between the 
parts $\not\equiv 0\pmod{15}$ is $\ge15$, with equality only if a part is 
$\equiv-2^3,-2^2,-2^1,-2^0\pmod{15}$, parts which are 
$\not\equiv -2^3,-2^2,-2^1,-2^0\pmod{15}$ are $>15$, the difference 
between the multiples of $15$ is $\ge60$, and the smallest multiple of $15$ is
\begin{equation}
\begin{cases}
\ge30+30\tau,&\text{ if $7$ is a part, }\\
\ge45+30\tau,&\text{ otherwise},
\end{cases}
\label{eq:1.1}
\end{equation}
where $\tau$ is number of non--multiples of $15$ in the partition.  Then
$$
G_3(n)=P_3(n).
$$
\label{thm:3}
\end{thm}
Theorem~\ref{thm:3} can be seen in \S3 as a special case of the following new remarkable 
four parameter {\it key identity}
\begin{align}
&\sum_{i,j,k,l}A^i B^j C^k D^l \sum_{\substack{i,j,k,l-\\ \text{constraints}}}
\frac{q^{T_\tau+T_{ab}+T_{ac}+\dots+T_{cd}-bc-bd-cd+4T_{Q-1}+3Q+2Q\tau}}
{(q)_a(q)_b(q)_c(q)_d(q)_{ab}(q)_{ac} (q)_{ad} (q)_{bc} (q)_{bd} (q)_{cd} (q)_Q}  \nonumber \\
&\cdot\left\{(1-q^a)+q^{a+bc+bd+Q}(1-q^b)+q^{a+bc+bd+Q+b+cd}\right\} \nonumber \\
&=(-Aq)_\infty(-Bq)_\infty(-Cq)_\infty(-Dq)_\infty
\label{eq:1.2}
\end{align}
under the transformations
\begin{equation}
\begin{cases}
\text{(dilation) }q\mapsto q^{15},\\
\text{(translations) }A\mapsto q^{-8},\ B\mapsto q^{-4},\ C\mapsto 
q^{-2},\ D\mapsto q^{-1}.
\label{eq:1.3}
\end{cases}
\end{equation}
{\bf Warning:} It is to be noted that in \eqn{1.2} and everywhere, $ab$, $ac$, $ad$, $bc$, 
$bd$, and $cd$ 
are parameters and that $ab$ is not $a$ multiplied by $b$, with similar interpretation 
for $ac$, $ad$, $bc$, $bd$ and $cd$.

In \eqn{1.2} we have made use of the standard notations
\begin{equation}
(a)_n=(a;q)_n=\begin{cases}
\prod^{n-1}_{j=0}(1-aq^j), &\text{ if } n>0, \\ 1, &\text{ if } n=0, \\ 
\prod^{-n}_{j=1}(1-aq^{-j})^{-1}, &\text{ if } n<0,
\end{cases}
\label{eq:1.4}
\end{equation}
 and
\begin{equation}
(a)_\infty=\lim_{n\to\infty}(a)_n,\text{ when } |q|<1.
\label{eq:1.5}
\end{equation}
In addition, in \eqn{1.2}  $T_n=\frac{n(n+1)}{2}$ is the $n$-th triangular number, and 
the $i,j,k,l-$constraints on the summation variables $a,b,c,d,ab,\dots,cd,Q$ are
\begin{equation}
\begin{cases}
i=a+ab+ac+ad+Q ,\\
j=b+ab+bc+bd+Q ,\\
k=c+ac+bc+cd+Q ,\\
l=d+ad+bd+cd+Q .
\end{cases}
\label{eq:1.6}
\end{equation}
The quantities $a,b,c,d,ab,\dots,cd$, and $Q$ may be interpreted as the number of 
parts occuring in certain colors as  will become 
clear in \S3. Finally in \eqn{1.2}, $\tau=a+b+c+d+ab+ac+ad+bc+bd+cd$.

Extracting the coefficients of $A^iB^jC^kD^l$  in \eqn{1.2}, we can rewrite  \eqn{1.2}
in the equivalent form as 
\begin{align}
&\sum_{\substack{i,j,k,l-\\ \text{constraints}}} 
\frac{q^{T_\tau+T_{ab}+T_{ac}+\dots+T_{cd}-bc-bd-cd+4T_{Q-1}+3Q+2Q\tau}}
{(q)_a(q)_b(q)_c(q)_d(q)_{ab}(q)_{ac} (q)_{ad} (q)_{bc} (q)_{bd} (q)_{cd} (q)_Q}  \nonumber \\
&\cdot\left\{(1-q^a)+q^{a+bc+bd+Q}(1-q^b)+q^{a+bc+bd+Q+b+cd}\right\} \nonumber \\
&=\frac{q^{T_i+T_j+T_k+T_l}} 
{(q)_i(q)_j(q)_k(q)_l},
\label{eq:1.7}
\end{align}
where we note that the sum in \eqn{1.7} is actually a finite $7-$fold sum.

In the next section we will first describe a two 
parameter refinement of Schur's theorem due to Alladi and Gordon \cite{16} 
(Theorem~\ref{thm:4} of \S2) who used the notion of partitions into colored integers. 
The analytic form of Theorem~\ref{thm:4} is the identity \eqn{2.7} of \S2.
Next in \S2 we will describe the generalization and three parameter 
refinement of G\"ollnitz's theorem due to Alladi, Andrews, and Gordon \cite{8} 
(Theorem~\ref{thm:5} of \S2), also obtained by using the notion of partitions into 
colored integers and by extending the method in \cite{16}. The analytic form of Theorem~\ref{thm:5} 
is the identity \eqn{2.15}. An advantage of this approach is that it 
shows clearly that the three parameter colored generalization of Theorem~\ref{thm:2}
is an extension of the two parameter colored generalization of Theorem~\ref{thm:1}. That is, 
Theorem~\ref{thm:4} is a special case of Theorem~\ref{thm:5}, and \eqn{2.7} a special 
case of \eqn{2.15} when any one of the parameters $i,j$, or $k$ is set equal to zero.

Pursuing the notion of reformulation in terms of partitions involving 
colored integers, we will construct in \S3 a four parameter colored generalization 
and refinement of Theorem~\ref{thm:3}, which is Theorem~\ref{thm:6}. In \S4 we will 
follow the development in \cite{12} and show that Theorem~\ref{thm:6} is the combinatorial
interpretation of the new identity \eqn{1.7}. The proof of this remarkable identity is given
in \S5 and \S6. If any one of the parameters $i,j,k$, or $l$ is set equal to $0$ in \eqn{1.7}, 
we get \eqn{2.15}; similarly, Theorem~\ref{thm:6} yields Theorem~\ref{thm:5} as a special case. 
Thus the question of Andrews is settled.

It is to be noted that in going from Theorem~\ref{thm:1} to Theorem~\ref{thm:2} or from 
Theorem~\ref{thm:4} to Theorem~\ref{thm:5}, the nature of the gap conditions does not change. 
Still the proof of Theorem~\ref{thm:5} is much deeper compared to that of Theorem~\ref{thm:4}. 
We wish to emphasize that in going from Theorem~\ref{thm:2} to Theorem~\ref{thm:3} or from 
Theorem~\ref{thm:5} to Theorem~\ref{thm:6}, the extension is non--routine because of two  very 
new conditions involving the multiplies of $15$ in Theorem~\ref{thm:3} and the parts in 
quaternary color in Theorem~\ref{thm:6}. More precisely, while the multiplies of $15$ 
(respectively parts in quaternary color)  satisfy special gap conditions among themselves, there is no
{\it interaction} between multiplies of $15$ (respectively parts in quaternary color) and the  the other
 parts  as far as gap conditions are concerned. Also, there is a subtle lower bound condition 
\eqn{1.1} in Theorem~\ref{thm:3} (respectively \eqn{3.3} in Theorem~\ref{thm:6}) on the multiplies of $15$ 
(respectively parts in quaternary color). These new conditions   might have been the reason that the problem 
posed by  Andrews remained unsolved for so long.

Alladi \cite{1}, \cite{5} and more recently Alladi and Berkovich \cite{13}, \cite{15} obtained 
several important consequences of \eqn{2.15} including new proofs and interpretations 
of Jacobi's celebrated triple product identity  as well as many new 
weighted partition identities. In a similar spirit, \eqn{1.2} leads to an extension of 
Jacobi's triple product identity which is stated as \eqn{7.2} in \S7. The new identity
\eqn{1.2} raises the exciting possibility of four parameter refinements of partition theorems 
of Capparelli \cite{23} and of Andrews--Bessenrodt--Olsson \cite{21}. We discuss this 
briefly in \S7 and plan to pursue this in detail later.

The new results (the identity \eqn{1.7}, Theorem~\ref{thm:3} and Theorem~\ref{thm:6}) 
were first announced without proofs in \cite{10}. The broad historical and mathematical significance 
of these results were discussed in the survey article \cite{Focus}.  Our main object here is to provide  the detailed proof of \eqn{1.7} and to show that 
\eqn{1.7} is the analytic form of Theorem~\ref{thm:6}.

\bigskip

\section{\bf Colored reformulation  and refinement \\of the Schur and G\"ollnitz theorems}

\medskip

We first describe the colored reinterpretation and refinement of Schur's theorem due 
to Alladi and Gordon \cite{16}.

For this purpose we assume that all positive integers occur in two primary 
colors $\Bbb A$ and $\Bbb B$ and that integers $\ge2$ occur also in the 
secondary color $\Bbb A\Bbb B$. Let $\Bbb A_n,\Bbb B_n$ and $\Bbb A\Bbb B_n$ 
denote the integer $n$ occuring in colors $\Bbb A,\Bbb B$, and $\Bbb 
A\Bbb B$, respectively.

In order to discuss partitions involving the colored integers, we need an 
ordering among them, and the one we choose is
\begin{equation}
\Bbb A_1<\Bbb B_1<\Bbb A\Bbb B_2<\Bbb A_2<\Bbb B_2<\Bbb A\Bbb B_3< \Bbb 
A_3<\Bbb B_3<\Bbb A\Bbb B_4<\cdots.
\label{eq:2.1}
\end{equation}
When the substitutions
\begin{equation}
\begin{cases}
\Bbb A_n\mapsto 3n-2,\ \Bbb B_n\mapsto 3n-1,\text{ for }n\ge1,\\
\Bbb A\Bbb B_n\mapsto 3n-3,\text{ for }n\ge2
\end{cases}
\label{eq:2.2}
\end{equation}
are made, the ordering \eqn{2.1} becomes
$$
1<2<3<4<\dots,
$$
the natural ordering among the positive integers.  Theorem~\ref{thm:4} stated below 
involving partitions into colored integers yields Theorem~\ref{thm:1} under the 
substitutions \eqn{2.2}.
Note that for a given integer $n$, the ordering is
\begin{equation}
\Bbb A\Bbb B_n<\Bbb A_n<\Bbb B_n.
\label{eq:2.3}
\end{equation}
In general, there are six orderings generated by the six permutations of the symbols 
$\Bbb A_n,\Bbb B_n$, and $\Bbb A\Bbb B_n $ as shown by Alladi--Gordon \cite{16},
but for the reasons given above, we have chosen \eqn{2.1}. 

In order to state Theorem~\ref{thm:4}, we define Type--$1$ partitions to be those of 
the form $\pi: m_1+m_2+\dots+m_\nu$ where the $m_i$ are colored integers 
from the list \eqn{2.1} such that $m_i-m_{i+1}\ge1$, where equality holds if either
\begin{equation}
\begin{cases}
m_i\text{ and } m_{i+1}\text{ are of the same primary color,}\\
\text{or $m_i$ is of a higher order color as given by \eqn{2.3}.}
\end{cases}
\label{eq:2.4}
\end{equation}
It is easy to check that under the substitutions \eqn{2.2}, the gap conditions in \eqn{2.4} 
translate to the difference conditions defining $G_1(n)$ in Theorem~\ref{thm:1}.

Next, for any partition $\pi$ into  colored integers, let $\nu_{\Bbb 
A}(\pi),  \nu_{\Bbb B}(\pi)$, and $\nu_{\Bbb A\Bbb B}(\pi)$ denote the 
number of parts of $\pi$ in colors $\Bbb A,\Bbb B$, and $\Bbb A\Bbb B$, 
respectively. We are now in a position to state the following result.
\begin{thm} (Alladi--Gordon \cite{16})\\

For given integers $i,j\ge0$, let $P_4(n;i,j)$ denote the number of 
partitions of $n$ into $i$ distinct parts all in color $\Bbb A$, and $j$ 
distinct parts all in color $\Bbb B$.

Let $G_4(n;a,b,ab)$ denote the number of Type--$1$ partitions $\pi$  of  $n$
with $a=\nu_{\Bbb A}(\pi), \ b=\nu_{\Bbb B}(\pi)$,  and  $ab=\nu_{\Bbb A\Bbb 
B}(\pi)$. Then
$$
\sum_{\substack{i=a+ab\\j=b+ab}} G_4(n;a,b,ab)=P_4(n;i,j).
$$
\label{thm:4}
\end{thm}
 It is clear that
\begin{equation}
\sum_{n} P_4(n;i,j)q^n= \frac{q^{T_i+T_j}} 
{(q)_i(q)_j}.
\label{eq:2.5}
\end{equation}
It turns out that (see \cite{16})
\begin{equation}
\sum_n G_4(n;a,b,ab)q^n=\frac{q^{T_{a+b+ab}+T_{ab}}}{(q)_a(q)_b(q)_{ab}}.
\label{eq:2.6}
\end{equation}
Thus Theorem~\ref{thm:4} is equivalent to the identity
\begin{equation}
\sum_{\substack{i=a+ab\\j= b+ab}} \frac{q^{T_{a+b+ab}+T_{ab}}}{(q)_a 
(q)_b(q)_{ab}}=\frac{q^{T_i+T_j}}{(q)_i(q)_j}.
\label{eq:2.7}
\end{equation}
A combinatorial proof of \eqn{2.7} was given  in \cite{16}. We remark that the identity \eqn{2.7} 
is a special case of $q$-Chu-Vandermonde summation formula (\cite{24}, (II.6), p. 236).

Multiplying \eqn{2.7} by $A^iB^j$ and summing the result over parameters $i$ and $j$ we can 
rewrite \eqn{2.7} in the equivalent form as 
\begin{equation}
\sum_{\substack{a,b,ab}}A^{a+ab}B^{b+ab} \frac{q^{T_{a+b+ab}+T_{ab}}}{(q)_a 
(q)_b(q)_{ab}}=(-Aq)_\infty(-Bq)_\infty.
\label{eq:ext.1}
\end{equation}
Note that the transformations
\begin{equation}
\begin{cases}
\text{(dilation) }q\mapsto q^3, \\ 
\text{(translations) }A\mapsto q^{-2},\ B\mapsto q^{-1},
\end{cases}
\label{eq:2.8}
\end{equation}
are equivalent to the substitutions \eqn{2.2} on the colored integers. Since under these 
substitutions, the gap conditions in \eqn{2.4} translate to the difference conditions  defining
$G_1(n)$ in Theorem~\ref{thm:1}, we see that under the transformations \eqn{2.8} the identity
\eqn{ext.1} becomes the analytic version of Theorem~\ref{thm:1}.
   
Pursuing the idea of partitions into colored integers, Alladi, Andrews, and 
Gordon \cite{8} reinterpreted and refined G\"ollnitz's theorem as we describe next.

Let us assume that the positive integers occur in three primary colors 
$\Bbb A,\Bbb B$, and $\Bbb C$, and that integers $\ge2$ occur also in the 
three secondary colors $\Bbb A\Bbb B, \Bbb A\Bbb C$, and $\Bbb B\Bbb C$. 
Each of the symbols $\Bbb A_n,\Bbb B_n,\Bbb C_n$, $\Bbb A\Bbb B_n$, $\Bbb 
A\Bbb C_n$, and $\Bbb B\Bbb C_n$, will represent the integer $n$ in the 
corresponding color.

As in the case of the colored reinterpretation of Schur's theorem, an ordering of 
these symbols is required here and the one we choose is
\begin{equation}
\Bbb A_1<\Bbb B_1<\Bbb C_1<\Bbb A\Bbb B_2<\Bbb A\Bbb C_2<\Bbb A_2<\Bbb 
B\Bbb C_2<\Bbb B_2<\Bbb C_2<\Bbb A\Bbb B_3<\Bbb A\Bbb C_3<\dots
\label{eq:2.9}
\end{equation}
One reason for this choice is that under the substitutions
\begin{equation}
\begin{cases}
\Bbb A_n\mapsto 6n-4,\ \Bbb B_n\mapsto 6n-2,\ \Bbb C_n\mapsto 6n-1, &\text{ 
for }n\ge1,\\
\Bbb A\Bbb B_n\mapsto 6n-6,\ \Bbb A\Bbb C_n\mapsto 6n-5,\ \Bbb B\Bbb 
C_n\mapsto 6n-3, &\text{ for } n\ge2,
\end{cases}
\label{eq:2.10}
\end{equation}
the ordering \eqn{2.9} becomes
$$
2<4<5<6<7<8<9<10<11<12<\dots,
$$
the natural ordering among the positive integers not equal to $1$ and $3$. 
Another reason is that Theorem~\ref{thm:5} stated below implies Theorem~\ref{thm:2} under the 
substitutions \eqn{2.10} as we shall soon see. We wish to comment that here too one can consider
other orderings of the colored integers and we refer the reader to \cite{8} for a discussion 
of the companion results to Theorem~\ref{thm:2} that these other orderings yield.

Note that for any given integer $n$, the ordering is
\begin{equation}
\Bbb A\Bbb B_n<\Bbb A\Bbb C_n<\Bbb A_n<\Bbb B\Bbb C_n<\Bbb B_n<\Bbb C_n.
\label{eq:2.11}
\end{equation}
In this case Type--$1$ partitions are defined as those of the form 
$\pi:m_1+m_2+\dots+m_\nu$, where the $m_i$ are colored integers from the 
list \eqn{2.9} such that $m_i-m_{i+1}\ge1$ with equality only if
\begin{equation}
\begin{cases}
m_i\text{ and } m_{i+1} \text{ are of the same primary color, or }\\
\text{ if $m_i$ is of a higher order color given by \eqn{2.11}}.
\end{cases}
\label{eq:2.12}
\end{equation}
We can now state the following:
\begin{thm} (Alladi--Andrews--Gordon \cite{8})\\

For given integers $i,j,k\ge0$, let $P_5(n,i,j,k)$ denote the number of 
partitions of $n$ into $i$ distinct parts in color $\Bbb A$, $j$ distinct 
parts in color $\Bbb B$, and $k$ distinct parts in color $\Bbb C$.

Let $G_5(n;a,b,c,ab,ac,bc)$ denote the number of Type--$1$ partitions $\pi$ of 
$n$ such that $\nu_{\Bbb A}(\pi)=a$, $\nu_{\Bbb B}(\pi)=b$, $\nu_{\Bbb C} 
(\pi)=c$, $\nu_{\Bbb A\Bbb B}(\pi)=ab$, $\nu_{\Bbb A\Bbb C}(\pi)=ac$ and 
$\nu_{\Bbb B\Bbb C}(\pi)=bc$. Then
$$
\sum_{\substack{i=a+ab+ac\\ j=b+ab+bc\\ k=c+ac+bc}} 
G_5(n;a,b,c,ab,ac,bc)=P_5(n;i,j,k). 
$$
\label{thm:5}
\end{thm}
Here too the notation involving $\nu$ and the parameters $ab,ac$, and $bc$ 
is as explained earlier in this section.
It is clear that
\begin{equation}
\sum_{n} P_5(n;i,j,k)q^n= 
\frac{q^{T_i+T_j+T_k}}{(q)_i(q)_j(q)_k}.
\label{eq:2.13}
\end{equation}
It is shown in \cite{8} that
\begin{equation}
\sum_n G_5(n,a,b,c,ab,ac,bc)q^n=\frac{q^{T_s+T_{ab}+T_{ac}+T_{bc-1}}(1-q^a+q^{a+bc})}
{(q)_a(q)_b(q)_c(q)_{ab}(q)_{ac}(q)_{bc}},
\label{eq:2.14}
\end{equation}
where $s=a+b+c+ab+ac+bc$ is the total number of parts in each Type--$1$ 
partition $\pi$. Thus Theorem~\ref{thm:5} is equivalent to the identity
\begin{equation}
\sum_{\substack{i=a+ab+ac\\ j=b+ab+bc\\ k=c+ac+bc}}
\frac{q^{T_s+T_{ab}+T_{ac}+T_{bc-1}}(1-q^a+q^{a+bc})} 
{(q)_a(q)_b(q)_c(q)_{ab}(q)_{ac}(q)_{bc}}= \frac{q^{T_i+T_j+T_k}} 
{(q)_i(q)_j(q)_k}.
\label{eq:2.15}
\end{equation}

The first proof of \eqn{2.15} in \cite{8} utilized  Whipple's $q$-analogue of Watson's 
transformation  ${}_8\varphi_7\rightarrow {}_4\varphi_3$ (\cite{24}, (III.18), p. 242),
the ${}_6\psi_6$ summation of Bailey (\cite{24}, (II.33), p. 239) as well as a transformation
formula ${}_3\varphi_2\rightarrow {}_3\varphi_2$ (\cite{24}, (III.9), p. 241).  
Subsequently in \cite{7}, Alladi and Andrews simplified the proof of \eqn{2.15}. This latter  
proof required only a special case of Jackson's $q$-analogue of Dougall's summation for 
${}_6\varphi_5$ (\cite{24}, (II.21), p. 238). In \cite{AR}, Riese used his computer algebra 
package qMultiSum to find a very simple recursive proof of \eqn{2.15}. In \cite{12}, Alladi and 
Berkovich proposed and proved a polynomial analog of \eqn{2.15}.

Multiplying \eqn{2.15} by $A^iB^jC^k$ and summing the result over parameters $i$, $j$ and $k$, we 
can rewrite \eqn{2.15} in the equivalent form as 
\begin{align}
\sum_{\substack{a,b,c,\\ab,ac,bc}}A^{a+ab+ac}B^{b+ab+bc}&C^{c+ac+bc} 
\frac{q^{T_s+T_{ab}+T_{ac}+T_{bc-1}}(1-q^a+q^{a+bc})} 
{(q)_a(q)_b(q)_c(q)_{ab}(q)_{ac}(q)_{bc}}\nonumber \\
&=(-Aq)_\infty(-Bq)_\infty(-Cq)_\infty,
\label{eq:ext.2}
\end{align}
with $s=a+b+c+ab+ac+bc$, as before. It is easy to verify that under the substitutions \eqn{2.10},
the gap conditions \eqn{2.12} become the difference conditions defining $G_2(n)$ in 
Theorem~\ref{thm:2}. Therefore, when transformations
\begin{equation}
\begin{cases}
\text{(dilation) } q\mapsto q^6,\\
\text{(translations) }A\mapsto q^{-4},\ B\mapsto q^{-2},\ C\mapsto q^{-1},
\end{cases}
\label{eq:2.16}
\end{equation}
are applied to the identity \eqn{ext.2}, it becomes the analytic version of Theorem~\ref{thm:2}.

When any one of the parameters $i,j$, or $k$ is set equal to $0$, \eqn{2.15} 
reduces \eqn{2.7}, and Theorem~\ref{thm:5} to Theorem~\ref{thm:4}. Pursuing this approach, we 
will describe in the next section how to construct a four parameter colored 
partition theorem, which reduces to Theorem~\ref{thm:5} when one of the parameters is set equal
to $0$ and will explain why this Theorem is equivalent to the new identity \eqn{1.7}.

It is to be noted that in Theorem~\ref{thm:5} the ternary color $\Bbb A\Bbb B\Bbb C$ 
is not utilized and so only a proper subset of the complete alphabet of $7$ 
 colors $\Bbb A,\Bbb B,\Bbb C,\Bbb A\Bbb B,\Bbb A\Bbb C,\Bbb B\Bbb 
C$, $\Bbb A\Bbb B\Bbb C$ is used. Although Andrews did not utilize the 
viewpoint of partitions into colored integers, starting with Schur's 
theorem he was able to construct \cite{17}, \cite{18} infinite hierarchies of partition 
theorems. His construction was based on choosing $r$ distinct residue classes $\pmod {2^r-1} $ 
and forming all possible sums of these residues to get the full set of  residue classes 
$\pmod {2^r-1} $. From the point of view presented here, the $r$ residue classes 
Andrews started with could be treated as $r$ primary colors, and the residue 
classes obtained by summation as secondary, ternary, $\dots$, colors 
depending on how many classes are summed. In other words, Andrews' construction is based on the
complete alphabet of colors. In contrast, our goal here is to construct 
a partition theorem starting with four primary colors $\Bbb A,\Bbb B,\Bbb 
C,\Bbb D$ and utilize only a proper subset of the alphabet of $15$ colors -- $4$ primary,
$6$ secondary, $4$ ternary, and $1$ quaternary.  It is due to the emphasis on 
selecting a proper subset of colors that one does not know how long the hierarchy 
beyond G\"ollnitz's Theorem extends. Indeed one did not know until recently whether 
such a partition theorem existed in case of the four primary colors, and if it 
did, whether it would reduce to Theorem~\ref{thm:5} when one of the colors is eliminated. 
The main realization in \cite{10} was that in addition to the primary and
secondary colors only the quaternary color $\Bbb A\Bbb B\Bbb C\Bbb D$ has to be retained, 
but all ternary colors $\Bbb A\Bbb B\Bbb C, \Bbb A\Bbb B\Bbb D,\Bbb A\Bbb 
C\Bbb D$, and $\Bbb B\Bbb C\Bbb D$ are to be discarded.  Without further 
ado we now describe the resolution of the Andrews problem along these lines.

\bigskip

\section{\bf A new four parameter partition theorem}

\medskip

We assume that all positive integers occur in four primary colors $\Bbb 
A,\Bbb B,\Bbb C$, and $\Bbb D$, that integers $\ge2$ also occur in the six 
secondary colors $\Bbb A\Bbb B,\Bbb A\Bbb C,\Bbb A\Bbb D,\Bbb B\Bbb C,\Bbb 
B\Bbb D$, and $\Bbb C\Bbb D$, and that integers $\ge4$ occur also in the 
quaternary color $\Bbb A\Bbb B\Bbb C\Bbb D$. As before the symbols $\Bbb 
A_n,\dots,\Bbb C \Bbb D_n,\Bbb A\Bbb B\Bbb C\Bbb D_n$ will represent the 
integer $n$ occuring in the corresponding colors.  We need now an ordering 
of the colored integers.

In the cases where we had $2$ or $3$ primary colors, Schur's theorem and G\"ollnitz's 
theorem helped us to select the orderings.  More precisely orderings \eqn{2.1} and \eqn{2.9} 
were chosen because under the substitutions \eqn{2.2} and \eqn{2.10} they yielded 
the natural ordering of the positive integers, and also the general 
Theorem~\ref{thm:4} and Theorem~\ref{thm:5} reduced to Theorem~\ref{thm:1} and 
Theorem~\ref{thm:2}, respectively.  In contrast, here in four dimensions, we do not (yet)
have a partition theorem to help us to select an ordering.  So what we do is to take the
ordering \eqn{2.11} and extend it as follows:
\begin{equation}
\begin{cases}
\text{if $m<n$ as ordinary (uncolored) positive integers},\\
\text{then $m$ in any color $<n$ in any color, and }\\
\text{if the same integer $n$ occurs in two different colors, the order is 
given by}\\
\Bbb A\Bbb B\Bbb C\Bbb D_n<\Bbb A\Bbb B_n<\Bbb A\Bbb C_n<\Bbb A\Bbb D_n<\Bbb A_n<\Bbb 
B\Bbb C_n<\Bbb B\Bbb D_n<\Bbb B_n<\Bbb C\Bbb D_n<\Bbb C_n<\Bbb D_n.
\end{cases}
\label{eq:3.1}
\end{equation}
With this ordering we define Type--$1$ partitions  as those of the form 
$\pi:m_1+m_2+\dots+m_\nu$, where the $m_i$ are colored integers ordered as in \eqn{3.1}, 
and such that parts in the non--quaternary colors differ by $\ge1$ 
with equality only if
\begin{equation}
\begin{cases}
\text{parts are of the same primary color, or}\\
\text{if the larger part is in a higher order color as given by 
\eqn{3.1}},
\end{cases}
\label{eq:3.2}
\end{equation}
and the gap between the parts in quaternary color is $\ge4$, with the added 
condition that the least quaternary part is
\begin{equation}
\begin{cases}
\ge3+2\tau, &\text{ if } \Bbb A_1\text{ is a part,}\\
\ge4+2\tau, &\text{ otherwise.}
\end{cases}
\label{eq:3.3}
\end{equation}
Here $\tau$ is the number of non--quaternary parts in the partition $\pi$.

The first important realization in our construction of Type--$1$ partitions 
is that all ternary parts are to be discarded.  The second major 
realization is that the quaternary parts do not {\it interact} with parts in 
the other colors in terms of gap conditions but only in terms  
of the lower bound \eqn{3.3}. For instance, for $n=12$ one may have the Type--$1$ partition 
${\Bbb A\Bbb B\Bbb C\Bbb D}_6+{\Bbb A}_6$, 
but not ${\Bbb A\Bbb B\Bbb C\Bbb D}_5+{\Bbb A}_7$. 
It is this special role played by the parts in quaternary 
color that causes a substantial increase in depth and difficulty in going 
beyond G\"ollnitz's theorem.  We are now in a position to state the new partition theorem.
\begin{thm} 

For given integers $i,j,k,l\ge0$, let $P_6(n;i,j,k,l)$ denote the number of 
partitions of $n$ into $i$ distinct parts in color $\Bbb A$, $j$ distinct 
parts in color $\Bbb B$, $k$ distinct parts in color $\Bbb C$, and $l$ 
distinct parts in color $\Bbb D$.

Let $\tilde{G}_6(n;a,b,c,d,ab,\cdots,cd,Q)$ denote the number of Type--$1$ partitions 
$\pi$ of $n$ such that $\nu_{\Bbb A}(\pi)=a,\dots, \nu_{\Bbb C\Bbb 
D}(\pi)=cd$, and $\nu_{\Bbb A\Bbb B\Bbb C\Bbb D}(\pi)=Q$, where $\nu_{\Bbb A}(\pi),
\\ \dots, \nu_{\Bbb C\Bbb D}(\pi), \nu_{\Bbb A\Bbb B\Bbb C\Bbb D}(\pi)$
denote the number of parts of $\pi$ in colors ${\Bbb A},\dots,{\Bbb C\Bbb D},
{\Bbb A\Bbb B\Bbb C\Bbb D}$, respectively. Then
$$
G_6(n;i,j,k,l)=P_6(n;i,j,k,l),
$$
where 
$$
G_6(n;i,j,k,l)=\sum_{\substack{i,j,k,l-\\ \text{constraints}}} 
\tilde{G}_6(n;a,b,c,d,ab,\dots,cd,Q)
$$
and the $i,j,k,l-$constraints on the parameters  $a, b,\dots, cd,Q$ are as in \eqn{1.6}.
\label{thm:6}
\end{thm}

For example, for $n=6$, $i=j=k=l=1$, there are ten pairs of partitions satisfying the
conditions in Theorem~\ref{thm:6}, as summarized in \tabref{1}.

\medskip

\begintable
$G_6(6;1,1,1,1)=10$ | $P_6(6;1,1,1,1)=10$ \crthick
$\Bbb A\Bbb B\Bbb C\Bbb D_6$       | $\Bbb A_1 + \Bbb B_1 + \Bbb C_1 + \Bbb D_3$ \cr
$\Bbb A\Bbb B_2 + \Bbb C\Bbb D_4$  | $\Bbb A_1 + \Bbb B_1 + \Bbb D_1 + \Bbb C_3$ \cr
$\Bbb A\Bbb C_2 + \Bbb B\Bbb D_4$  | $\Bbb A_1 + \Bbb C_1 + \Bbb D_1 + \Bbb B_3$ \cr
$\Bbb A\Bbb D_2 + \Bbb B\Bbb C_4$  | $\Bbb B_1 + \Bbb C_1 + \Bbb D_1 + \Bbb A_3$ \cr
$\Bbb B\Bbb C_2 + \Bbb A\Bbb D_4$  | $\Bbb A_1 + \Bbb B_1 + \Bbb C_2 + \Bbb D_2$ \cr
$\Bbb B\Bbb D_2 + \Bbb A\Bbb C_4$  | $\Bbb A_1 + \Bbb C_1 + \Bbb B_2 + \Bbb D_2$ \cr
$\Bbb C\Bbb D_2 + \Bbb A\Bbb B_4$  | $\Bbb A_1 + \Bbb D_1 + \Bbb B_2 + \Bbb C_2$ \cr
$\Bbb A_1  + \Bbb B_2 + \Bbb C\Bbb D_3$|$\Bbb B_1 + \Bbb C_1 + \Bbb A_2 + \Bbb D_2$\cr
$\Bbb A_1  + \Bbb B\Bbb C_2 + \Bbb D_3$|$\Bbb B_1 + \Bbb D_1 + \Bbb A_2 + \Bbb C_2$\cr
$\Bbb A_1  + \Bbb B\Bbb D_2 + \Bbb C_3$|$\Bbb C_1 + \Bbb D_1 + \Bbb A_2 + \Bbb B_2$ 
\endtable
\mycap{1}{}

A strong four parameter refinement of Theorem~\ref{thm:3} follows from Theorem~\ref{thm:6}  
under the substitutions
\begin{equation}
\begin{cases}
\Bbb A_n\mapsto 15n-8,\Bbb B_n\mapsto 15n-4,\Bbb C_n\mapsto 15n-2,\Bbb 
D_n\mapsto 15n-1,\text{ for } n\ge1,\\
\text{and consequently } \Bbb A\Bbb B_n\mapsto 15n-12, \Bbb A\Bbb 
C_n\mapsto 15n-10, \Bbb A \Bbb D_n\mapsto 15n-9,\\
\Bbb B\Bbb C_n\mapsto 15n-6,\Bbb B\Bbb D_n\mapsto 15n-5,\Bbb C\Bbb D 
\mapsto 15n-3,\text{ for } n\ge2,\\
\text{and }\Bbb A\Bbb B\Bbb C\Bbb D_n\mapsto 15n-15,\text{ for } n\ge4.
\end{cases}
\label{eq:3.4}
\end{equation}
The nice thing about these substitutions is that they convert the ordering \eqn{3.1} to
$$
7<11<13<14<18<20<21<22<24<25<26<27<28<29<\dots,
$$
the natural ordering  among  integers  $\not\equiv 2^3,2^2,2^1,2^0\pmod{15}$  and 
$\ne 3,5,6,9,10,
\\12,15,30 $. The substitutions \eqn{3.4} imply that the primary colors
correspond to the residue classes
$$
-2^3,-2^2,-2^1,-2^0\pmod{15}.
$$
Since
$$
2^3+2^2+2^1+2^0=15,
$$
the ternary colors correspond to the residue classes
$$
2^3,2^2,2^1,2^0\pmod{15},
$$
which are precisely the four residue classes not considered in Theorem~\ref{thm:3}. 
Also since the residue classes relatively prime to $15$ are 
$\pm2^3,\pm2^2,\pm2^1,\pm 2^0$ $\pmod{15}$, it follows that the integers in 
secondary colors correspond to the non--multiples of $15$ that are not 
relatively prime to $15$. Finally the integers in quaternary color correspond 
to the multiples of $15$, which are $\ge45$.  These features make Theorem~\ref{thm:3} 
especially appealing.

In the next section we show that Theorem~\ref{thm:6} is combinatorially equivalent to 
\eqn{1.7}.

\bigskip

\section{\bf Combinatorial interpretation of the identity \eqn{1.7}}

\medskip
In this section we will slightly generalize the analysis in \S4 of \cite{12}
to show that the identity \eqn{1.7} is an analytic version of  Theorem~\ref{thm:6}.  
First we observe that if the colors $\Bbb A,\Bbb B,\Bbb C,\Bbb D$ are ordered in some fashion 
such as in \eqn{3.1}, then it is clear that
\begin{equation}
\sum_nP_6(n;i,j,k,l)q^n=\frac{q^{T_i+T_j+T_k+T_l}}{(q)_i(q)_j(q)_k(q)_l},
\label{eq:4.1}
\end{equation}
because for any integer $i\ge0$,
$$
\frac{q^{T_i}}{(q)_i}
$$
is the generating function of partitions into $i$ distinct positive parts.  We 
will now show that the generating function of Type--$1$ partitions $\pi$ with 
$\nu_{\Bbb A}(\pi)=a,\dots,\nu_{\Bbb C\Bbb D}(\pi)\\
=cd$, and $\nu_{\Bbb A\Bbb B\Bbb C\Bbb D}(\pi)=Q$ is the  summand on the left of the 
identity \eqn{1.7}. To this end we first observe that for $i\ge0$,
$$
\frac{q^{T_i-i}}{(q)_i}
$$
is the generating function of partitions into $i$ distinct nonnegative parts and 
$$
\frac1{(q)_i}, \text{   }\frac{q^i}{(q)_i},\text{   }\frac1{(q)_{i-1}}=\frac{(1-q^i)}{(q)_i}
$$
are the generating functions of partitions into $i$ nonnegative parts with the least part 
$\ge 0$, the least part $>0$, the least part $=0$, respectively. 
Since the discussion of the least part will play a 
prominent role, we introduce the notation $\lambda(\pi)$ for the least part 
of the partition $\pi$.

The discussion of Type--$1$ partitions $\pi$ splits into three cases.

\noindent
\underbar{Case 1}: $\lambda(\pi)=\Bbb A_1$

Decompose $\pi$ as $\tilde\pi\cup\pi'$, where $\tilde\pi$ has all the parts 
of $\pi$ in primary and secondary colors, and $\pi'$ has the remaining 
parts all in quaternary color.  Note that 
$\nu(\tilde\pi)=\tau=a+b+c+d+ab+\dots+cd$, and $\nu(\pi')=Q$.

Subtract $1$ from the smallest part of $\tilde\pi$, $2$ from the second 
smallest part of $\tilde\pi,\dots$, and $\tau$ from the largest part of 
$\tilde\pi$. We call this process the {\it Euler subtraction}. It accounts 
for the term $q^{T_\tau}$ on the left side of \eqn{1.7}. Note that Euler 
subtraction preserves the ordering among the colored parts and therefore can be
easily reversed. After this  subtraction, $\tilde\pi$ can be decomposed into ten 
monochromatic partitions $\tilde\pi_{\Bbb A},\tilde\pi_{\Bbb B},\dots,\tilde\pi_{\Bbb C\Bbb D}$, 
with the color indicated by the subscript, and satisfying the following 
conditions:
\begin{equation}
\begin{cases}
\lambda(\tilde\pi_{\Bbb A})=0, & \nu(\tilde\pi_{\Bbb A})=a,\\
\lambda(\tilde\pi_{\Bbb B})\ge0, & \nu(\tilde\pi_{\Bbb B})=b,\\
\lambda(\tilde\pi_{\Bbb C})\ge0, & \nu(\tilde\pi_{\Bbb C})=c,\\
\lambda(\tilde\pi_{\Bbb D})\ge0, & \nu(\tilde\pi_{\Bbb D})=d,
\end{cases}
\label{eq:4.2}
\end{equation}
\begin{equation}
\begin{cases}
\tilde\pi_{\Bbb A\Bbb B}\text{ has $ab$ distinct (positive) parts,}\\
\tilde\pi_{\Bbb A\Bbb C}\text{ has $ac$ distinct (positive) parts,}\\
\tilde\pi_{\Bbb A\Bbb D}\text{ has $ad$ distinct (positive) parts,}
\end{cases}
\label{eq:4.3}
\end{equation}
\begin{equation}
\begin{cases}
\tilde\pi_{\Bbb B\Bbb C}\text{ has $bc$ distinct non--negative parts,}\\
\tilde\pi_{\Bbb B\Bbb D}\text{ has $bd$ distinct non--negative parts,}\\
\tilde\pi_{\Bbb C\Bbb D}\text{ has $cd$ distinct non--negative 
parts.}
\end{cases}
\label{eq:4.4}
\end{equation}
Note that in a Type--$1$ partition the gap between the parts in the 
same secondary color is $\ge2$, therefore after the Euler subtraction each 
monochromatic partition with the parts in secondary color will have distinct 
parts.  However, since the secondary colors $\Bbb B\Bbb C,\Bbb B\Bbb D$, 
and $\Bbb C\Bbb D$ are of higher order than the primary color $\Bbb A$ (see 
\eqn{3.1}), the partitions $\tilde\pi_{\Bbb B\Bbb C},\tilde\pi_{\Bbb B\Bbb D}$ and 
$\tilde\pi_{\Bbb C\Bbb D}$, could have least part equal to $0$. For example, if from 
the Type--$1$ partition $\Bbb A_1\Bbb B\Bbb C_2$ we subtract $1$ from $\Bbb 
A_1$ and $2$ from $\Bbb B\Bbb C_2$, we get $\Bbb A_0\Bbb B\Bbb C_0$. Such a 
possibility does not arise with the colors $\Bbb A\Bbb B,\Bbb A\Bbb C$, and 
$\Bbb A\Bbb D$ which are of lower order than $\Bbb A$. Hence the conditions 
\eqn{4.3} are slightly different from the conditions \eqn{4.4}.

Our most crucial observation concerning these monochromatic partitions is 
that they are {\it independent} of each other. 
Thus we may calculate the generating function of each of these 
monochromatic partitions.
The conditions \eqn{4.2} imply that the generating functions of $\tilde\pi_{\Bbb 
A}, \tilde\pi_{\Bbb B}, \tilde\pi_{\Bbb C}$, and $\tilde\pi_{\Bbb D}$ are
\begin{equation}
\frac{(1-q^a)}{(q)_a},\frac1{(q)_b},\frac1{(q)_c}, \frac1{(q)_d},
\label{eq:4.5}
\end{equation}
 respectively.
Similarly, \eqn{4.3} implies that the generating functions of $\tilde\pi_{\Bbb 
A\Bbb B}$, $\tilde\pi_{\Bbb A\Bbb C}$, and $\tilde\pi_{\Bbb A\Bbb D}$ are
\begin{equation}
\frac{q^{T_{ab}}}{(q)_{ab}}, \frac{q^{T_{ac}}}{(q)_{ac}} ,\frac{q^{T_{ad}}}{(q)_{ad}},
\label{eq:4.6}
\end{equation}
respectively.
Finally, \eqn{4.4} implies that the generating functions of $\tilde\pi_{\Bbb 
B\Bbb C},\tilde\pi_{\Bbb B\Bbb D}$, and $\tilde\pi_{\Bbb C\Bbb D}$ are
\begin{equation}
\frac{q^{T_{bc}-bc}}{(q)_{bc}},\frac{q^{T_{bd}-bd}}{(q)_{bd}}, \frac{q^{T_{cd}-cd}}{(q)_{cd}},
\label{eq:4.7}
\end{equation}
respectively.
The parts in quaternary color have no interaction with the other parts in 
primary and secondary colors and so act {\it independently} except for the lower bound,
which for  $\lambda(\pi)=\Bbb A_1$ is  
$$
\lambda(\pi')\ge3+2\tau .
$$
 Since the gap between the parts in 
 quaternary color is $\ge4$, the generating function of $\pi'$ is given by 
\begin{equation}
\frac{q^{4T_{Q-1}+3Q+2Q\tau}}{(q)_Q}.
\label{eq:4.8}
\end{equation}
Taking the product of the terms in \eqn{4.5}--\eqn{4.8} along with $q^{T_\tau}$ (to 
account for the Euler subtraction) we see that if $a,b,c,d,ab,\dots,cd$, 
and $Q$ are specified, then the generating function of Type--$1$ partitions with 
$\lambda(\pi)=\Bbb A_1$ is
\begin{equation}
\frac{q^{T_\tau+T_{ab}+\dots+T_{cd}-bc-bd-cd}(1-q^a)q^{4T_{Q-1}+3Q+ 
2Q\tau}}{(q)_a(q)_b(q)_c(q)_d(q)_{ab}(q)_{ac} (q)_{ad} (q)_{bc} (q)_{bd} (q)_{cd} (q)_Q}.  
\label{eq:4.9}
\end{equation}

\noindent
\underbar{Case 2}: $\lambda(\pi)=\Bbb B_1$.

Here the discussion is identical to Case $1$ except that after the Euler subtraction
the smallest permissible value of parts in colors $\Bbb A,\Bbb B\Bbb C$, and $\Bbb B\Bbb 
D$, will be $1$ because these colors are of lower order compared to $\Bbb B$  (but they 
are not of lower order than $\Bbb A$). 
Also by the second case of \eqn{3.3} the smallest permissible value of a part in quaternary 
color is now $4+2\tau$ which is one higher than in Case $1$.  Finally, 
the term
\begin{equation}
\frac1{(q)_{b-1}}=\frac{(1-q^b)}{(q)_b}
\label{eq:4.10}
\end{equation}
will now replace $\frac1{(q)_{a-1}}$, because $\lambda(\tilde\pi_{\Bbb B})=0$ and
$\lambda(\tilde\pi_{\Bbb A})\ge0$. Thus if $a,b,c,d,ab,\dots,
\\cd$ and $Q$ are specified, then from \eqn{4.9}, \eqn{4.10} and the preceding observations 
we see that the generating function of Type--$1$ partitions with $\lambda(\pi)=\Bbb B_1$ is
\begin{equation}
\frac{q^{T_\tau+T_{ab}+\dots+T_{cd}-bc-bd-cd}(1-q^b)q^{4T_{Q-1}+3Q+ 
2Q\tau} 
q^{a+bc+bd+Q}}{(q)_a(q)_b(q)_c(q)_d(q)_{ab}(q)_{ac} (q)_{ad} (q)_{bc} (q)_{bd} (q)_{cd} (q)_Q}   .
\label{eq:4.11}
\end{equation}

\noindent
\underbar{Case 3}: $\lambda(\pi)>\Bbb B_1$ (all cases other than Cases $1$ and $2$)

Note that although integers $n\ge2$ can occur in colors other than the 
primary color, the integer $1$ can occur only in the primary colors $\Bbb 
A<\Bbb B<\Bbb C<\Bbb D$. Thus for a Type--$1$ partition, the smallest part 
can be either $\Bbb A_1$ (which is Case $1$), or $\Bbb B_1$ (which is Case 
$2$), or $>\Bbb B_1$ (which is Case $3$). In this final case, the discussion is 
the same as in Case $1$ except that after the Euler subtraction
the smallest permissible value of parts in colors $\Bbb A,\Bbb B\Bbb C,\Bbb B\Bbb D
,\Bbb C\Bbb D$, and $Q$ will be  $1$ .  Thus if $a,b,c,d,ab,\dots,cd$, and $Q$ are 
specified, then the generating function of Type--$1$ partitions with 
$\lambda(\pi)>\Bbb B_1$ is
\begin{equation}
\frac{q^{T_\tau+T_{ab}+\dots+T_{cd}-bc-bd-cd} q^{4T_{Q-1}+3Q+ 
2Q\tau} 
q^{a+bc+bd+Q+b+cd}}{(q)_a(q)_b(q)_c(q)_d(q)_{ab}(q)_{ac} (q)_{ad} (q)_{bc} (q)_{bd} (q)_{cd} (q)_Q}  .
\label{eq:4.12}
\end{equation}

Finally, adding the expressions in \eqn{4.9}, \eqn{4.11}, and \eqn{4.12}
and summing the result  over parameters $a,b,c,d,\dots,cd$, and $Q$ satisfying the $i,j,k,l-$
constraints \eqn{1.6} yields the  sum on the left in \eqn{1.7}. This is precisely 
the generating function of the sum over Type--$1$ partitions in Theorem~\ref{thm:6}.
Comparing this with \eqn{4.1} we see that Theorem~\ref{thm:6} is the combinatorial
version of the identity \eqn{1.7}.

Thus to prove Theorem~\ref{thm:6}, we need to establish \eqn{1.2}, which is equivalent to 
\eqn{1.7}. This we will  accomplish in the next two sections.

\bigskip

\section{\bf The key identity \eqn{1.2} as a constant term identity}

\medskip

In this section we will follow the development in \S4 of \cite{8} and 
transform \eqn{1.2} into the equivalent constant term 
identity \eqn{5.18} below.  To state this identity and for later use, we need 
the following notations and definitions:
\begin{equation}
[z^m]P(z)=\text{coefficient of $z^m$ in the Laurent expansion of $P(z)$},
\label{eq:5.1}
\end{equation}
\begin{equation}
(A_1,A_2,\dots,A_r; q)_n=(A_1,A_2,\dots,A_r)_n=(A_1)_n(A_2)_n\dots 
(A_r)_n,
\label{eq:5.3}
\end{equation}
\begin{equation}
{}_r\psi_r\left(\begin{matrix} a_1,a_2,\dots,a_r\\ b_1, b_2,\dots,b_r\end{matrix}; 
q,z\right)=\sum^\infty_{n=-\infty} 
\frac{(a_1,a_2,\dots,a_r)_n}{(b_1,b_2,\dots,b_r)_n}z^n,
\label{eq:5.4}
\end{equation}
\begin{equation}
{}_{r+1}\phi_r\left(\begin{matrix} a_1,a_2,\dots,a_{r+1}\\ b_1, b_2,\dots,b_r\end{matrix}; 
q,z\right)=\sum^\infty_{n=0} 
\frac{(a_1,a_2,\dots,a_{r+1})_n}{(q,b_1,\dots,b_r)_n}z^n.
\label{eq:qhyp}
\end{equation}
A ${}_{r+1}\phi_r$ series is called balanced if $z=q$ and $a_1\ldots a_{r+1}q=b_1\ldots b_r$.
We will also need very-well-poised series defined as
\begin{equation}
{}_{r+1}W_r(a_1;a_4,\dots,a_{r+1};q,z)=
{}_{r+1}\phi_r\left(\begin{matrix} a_1,q\sqrt{a_1},-q\sqrt{a_1},a_4,\dots,a_{r+1}
\\ \sqrt{a_1},-\sqrt{a_1},\frac{qa_1}{a_4},\dots,\frac{qa_1}{a_{r+1}}\end{matrix}; 
q,z\right).
\label{eq:WP}
\end{equation}

Next, we will make use of the relation among triangular numbers
\begin{equation}
T_\tau+4T_{Q-1}+(3+2\tau)Q=T_{\tau+2Q}
\label{eq:5.5}
\end{equation}
and get rid of the condition that
\begin{equation}
\tau=a+b+c+d+ab+ac+ad+bc+bd+cd 
\label{eq:5.6}
\end{equation}
in \eqn{1.2} by rewriting it in the form
$$
[z^0]\sum_{\substack{a,b,\ldots,cd,Q\ge 0 \\ -\infty < t < \infty}} 
\frac{q^{T_t+T_{ab}+\dots +T_{cd}-bc-bd-cd}} {(q)_a(q)_b\dots(q)_{cd}(q)_Q}
$$
\begin{equation}
\cdot\big\{1-q^a+q^{a+bc+bd+Q}(1-q^b)+q^{a+bc+bd+Q+b+cd}\big\}
\label{eq:5.7}
\end{equation}
$$
\cdot z^{(a+b+\dots+cd+2Q-t)} A^{a+ab+ac+ad+Q}\dots D^{d+ad+bd+cd+Q}
=(-Aq,-Bq,-Cq,-Dq)_\infty.
$$

At this point we observe that all twelve sums in \eqn{5.7} over the variables 
$t,a,b,\dots,
\\Q$ can be evaluated to be infinite products. To deal with 
the sum over $t$, we use Jacobi's triple product identity (\cite{24}, (II.28), p. 239)
\begin{equation}
\sum^\infty_{t=-\infty} q^{T_t}z^{-t}=(q,-z,-\frac qz)_\infty.
\label{eq:5.8}
\end{equation}
To deal with the sums over $a,b,c,d,Q$ and $ab,ac,ad,bc,bd,cd,$ we use two 
identities of Euler
\begin{equation}
\sum_{n\ge0}\frac{z^n}{(q)_n}=\frac1{(z)_\infty},
\label{eq:5.9}
\end{equation}
and
\begin{equation}
\sum_{n\ge0}\frac{q^{T_n}}{(q)_n}z^n=(-qz)_\infty,
\label{eq:5.10}
\end{equation}
respectively.  This way we derive
\begin{align}
&[z^0]\left\{Az \frac{(q,-z,-\frac qz,-ABqz,-ACqz,-ADqz,-BCz,-BDz, -CDz)_\infty}
{(Az,Bz,Cz,Dz,ABCDz^2,-Aq,-Bq,-Cq,-Dq)_\infty}\right.+\nonumber\\
&\left.(1+BCDz^2) \frac{(q,-z,-\frac qz,-ABqz,-ACqz,-ADqz,-BCqz,-BDqz,-CDqz)_\infty} {(Aqz,Bz,Cz,Dz, 
ABCDqz^2,-Aq,-Bq,-Cq,-Dq)_\infty}\right\} \nonumber\\ &=1.
\label{eq:5.11}
\end{align}

To proceed further we recall the famous Bailey ${}_6\psi_6$ summation 
(\cite{24}, (II.33), p. 239)
\begin{align}
{}_6\psi_6&={}_6\psi_6\left(\begin{matrix}
q\sqrt{-z},-q\sqrt{-z},-A^{-1},-B^{-1},-C^{-1},-D^{-1}\\
\sqrt{-z},-\sqrt{-z},Aqz,Bqz,Cqz,Dqz \end{matrix};q,ABCDqz^2\right)\label{eq:5.12}\\
&=\frac{(q,-qz,-\frac qz,-ABqz,-ACqz,-ADqz,-BCqz,-BDqz,-CDqz)_\infty} 
{(Aqz,Bqz,Cqz,Dqz,ABCDqz^2,-Aq,-Bq,-Cq,-Dq)_\infty},\nonumber
\end{align}
provided $|qz^2ABCD|<1$. Thanks to \eqn{5.12} we can rewrite \eqn{5.11} as
\begin{equation}
[z^0]\{R(A,B,C,D,z) \frac{(1+z)}{(1-Az)(1-Bz)(1-Cz)(1-Dz)} 
{}_6\psi_6\}=1,
\label{eq:5.13}
\end{equation}
where
$$
R(A,B,C,D,z)=(1-Az)(1+BCDz^2)
$$
\begin{equation}
+\frac{Az}{1-ABCDz^2}(1+BCz)(1+BDz)(1+CDz).
\label{eq:5.14}
\end{equation}

Next, we use the easily verifiable relation
\begin{equation}
(a)_{-n}= 
\frac{(-1)^na^{-n}q^{T_n}}{\left(\frac qa\right)_n},n>0
\label{eq:5.15}
\end{equation}
to decompose the function ${}_6\psi_6$ in \eqn{5.13} as
\begin{align}
&{}_6\psi_6=\sum_{n\ge1} \frac{1+zq^{2n}}{1+z} 
\frac{(-A^{-1},-B^{-1},-C^{-1},-D^{-1})_n}{(Aqz,Bqz,Cqz,Dqz)_n} (ABCDqz^2)^n+1\nonumber\\
+&\frac{(1-Az)(1-Bz)(1-Cz)(1-Dz)}{1+z}\sum_{n\ge1} (z+q^{2n}) 
\frac{\left(\frac q{Az},\frac q{Bz},\frac q{Cz},\frac q{Dz}\right)_{n-1}}
{(-Aq,-Bq,-Cq,-Dq)_n}\nonumber\\
&\qquad\qquad\qquad \cdot (ABCD)^{n-1}q^nz^{2n-4}.
\label{eq:5.16}
\end{align}
Note that the first sum on the right of \eqn{5.16} contains only positive powers 
of $z$ and therefore, can be discarded in \eqn{5.13}.
Since
\begin{equation}
[z^0]\left\{R(A,B,C,D,z) \frac{(1+z)}{(1-Az)(1-Bz)(1-Cz)(1-Dz)}\right\}=1 ,
\label{eq:5.17}
\end{equation}
we derive from \eqn{5.13} the desired {\it constant term identity}
\begin{align}
[z^0]\bigg\{R(A,B,C,D,z)\sum_{n\ge1}&(z+q^{2n})(ABCD)^{n-1}q^nz^{2n-4}\nonumber\\
&\cdot\frac{\left(\frac q{Az},\frac q{Bz},\frac q{Cz},\frac q{Dz}\right)_{n-1}}
{(-Aq,-Bq,-Cq,-Dq)_n}\bigg\}=0,
\label{eq:5.18}
\end{align}
where $R(A,B,C,D,z)$ is defined in \eqn{5.14}.
Notice that the sum on the left of \eqn{5.18} can be rewritten as
\begin{equation}
\frac{q(z+q^2)}{z^2(1+Aq)(1+Bq)(1+Cq)(1+Dq)}{}_8W_7(-\frac{q^2}{z};
q,\frac{q}{Az},\frac{q}{Bz},\frac{q}{Cz},\frac{q}{Dz};q,qz^2F),
\label{eq:ins1}
\end{equation}
where $F=ABCD$.
If we use the transformation formula (III.24) in \cite{24} together with \eqn{5.18} and 
\eqn{ins1} we obtain
$$
[z^0]\left\{\frac{R(A,B,C,D,z)(-\frac{q^2}{z},q^2zABC,q^2zABD,q^2zACD,q^2zBCD)_\infty}
{z(-q^3zF,qz^2F)_\infty}\right.
$$
\begin{equation}
{}_8W_7(-q^2zF;-qA,-qB,-qC,-qD,qz^2F;q,q)\bigg\}=0.
\label{eq:ins2}
\end{equation}
It is interesting to observe that ${}_8W_7(-q^2zF;-qA,-qB,-qC,-qD,qz^2F;q,q)$ in \eqn{ins2} 
is very-well-poised and balanced at the same time.

In spite of all our attempts to prove \eqn{5.18} or \eqn{ins2} in a purely \qhg 
fashion, we succeeded only when one of the parameters, say $D$, 
is set to be zero. In the latter case, we obtain from \eqn{ins2} with $D=0$ the following
\begin{equation}
[z^0]H(A,B,C,z)=0,
\label{eq:ins3}
\end{equation}
where
\begin{equation}
H(A,B,C,z)=(\frac{1}{z}+ABCz)(-\frac{q^2}{z},q^2zABC)_\infty
{}_3\phi_2\left(\begin{matrix} -qA,-qB,-qC\\ q^2zABC, -\frac{q^2}{z}\end{matrix};q,q\right). 
\label{eq:ins4}
\end{equation}
To prove \eqn{ins3} we follow \cite{8} and notice that
\begin{equation}
H(A,B,C,z)=-H(A,B,C,-\frac{1}{ABCz}).
\label{eq:ins5}
\end{equation}

Our strategy to prove \eqn{5.18} in full generality is briefly
described below.

First, we multiply both sides of \eqn{5.18} by $(-Aq,-Bq,-Cq,-Dq)_\infty$ to obtain
$$
[z^0]\{R(A,B,C,D,z)\sum_{n\ge1}(z+q^{2n})F^{n-1}q^nz^{2n-4}
$$
\begin{equation}
\cdot\left(\frac q{Az},\frac q{Bz},\frac q{Cz},\frac q{Dz}\right)_{n-1} 
(-Aq^{n+1},-Bq^{n+1},-Cq^{n+1},-Dq^{n+1})_\infty\}=0.
\label{eq:5.19}
\end{equation}
Next we rewrite \eqn{5.19} as
\begin{equation}
\lim_{m\to\infty}S(m)=0,
\label{eq:5.20}
\end{equation}
where
\begin{equation}
S(m)=[z^0]\{R(A,B,C,D,z)\sum^m_{n=1}(z+q^{2n})F^{n-1}q^n 
\label{eq:5.21}
\end{equation}
$$
\cdot z^{2n-4}\left(\frac q{Az},\frac q{Bz},\frac q{Cz},\frac q{Dz}\right)_{n-1} 
\left(-Aq^{n+1},-Bq^{n+1},-Cq^{n+1},-Dq^{n+1}\right)_{m-n}\}. 
$$
In the next section we will establish the identity \eqn{6.1} for $S(m)$, which 
would enable us to prove \eqn{5.20} and, as a result, \eqn{1.2}.

\bigskip

\section{\bf Proof of the constant term identity \eqn{5.18}}

\medskip

In this section we will prove that
\begin{align}
S(m)&=\sum_{2m+1\le i+j+k+l} q^{T_i+T_j+T_k+T_l} A^i B^jC^kD^l
\qBin{m}{i} \qBin{m}{j} \qBin{m}{k} \qBin{m}{l} \label{eq:6.1} \\
&-\sum_{L\ge1} q^{2L(m+1)}F^L\sum_{i+j+k+l=2m-2L} 
q^{T_i+T_j+T_k+T_l} A^iB^jC^k D^l
\qBin{m}{i} \qBin{m}{j} \qBin{m}{k} \qBin{m}{l},\nonumber
\end{align}
where $F=ABCD$, $S(m)$ is defined in \eqn{5.21}, and $q-$binomial coefficients are defined by
\begin{equation}
\qBin{m}{i}=
\begin{cases}
\frac{(q)_m}{(q)_i(q)_{m-i}},&\text{ if }0 \le i\le m,\\
0,&\text{ otherwise}.
\end{cases}
\label{eq:6.3}
\end{equation}
This is accomplished by showing that both sides in \eqn{6.1} satisfy the same 
recurrence relations and the same initial conditions. 

Before we move on, we 
remark that setting one of the parameters $A,B,C,D$ to zero will turn the 
second sum on the right of \eqn{6.1} into zero.  In other words, the reduction 
from the four parameter case to the three parameter (G\"ollnitz) case is rather 
dramatic.

Next, we note that while it is not immediately obvious from \eqn{5.14}, 
the rational function $R(A,B,C,D,z)$ is actually symmetric in $A,B,C,D$. 
To render this symmetry explicit we  expand the term $(1-ABCD z^2)^{-1}$ 
in \eqn{5.14} to obtain after rearrangement
$$
R(A,B,C,D,z)=1-F z^3+(A+B+C+D)\sum_{L\ge1}F^L z^{2L+1}
$$
\begin{equation} 
+(ABC+ABD+ACD+BCD)\sum_{L\ge0} F^L z^{2L+2}.
\label{eq:6.4}
\end{equation}

Next from the definition of $S(m)$ in \eqn{5.21} it follows that
\begin{equation}
S(m)-(1+Aq^m)(1+Bq^m)(1+Cq^m)(1+Dq^m)S(m-1) 
\label{eq:6.5}
\end{equation}
$$
=[z^0]\bigg\{R(A,B,C,D,z)(z+q^{2m})F^{m-1}q^m z^{2m-4}
\left(\frac q{Az},\frac q{Bz},\frac q{Cz},\frac q{Dz}\right)_{m-1}\bigg\},
$$
because in this subtraction only the term corresponding to $n=m$ in the sum 
\eqn{5.21} survives.  At this stage we expand the product
$$
\left(\frac q{Az},\frac q{Bz},\frac q{Cz},\frac q{Dz}\right)_{m-1}
$$
using the $q-$binominal theorem
\begin{equation} 
(zq)_N=\sum^N_{\nu=0}(-1)^\nu z^\nu q^{T_\nu}\qBin{N}{\nu}
=\sum^N_{\nu=0}(-1)^{N-\nu} z^{N-\nu} q^{T_{N-\nu}} \qBin{N}{\nu}
\label{eq:6.6}
\end{equation}
to get
$$
F^{m-1}q^m z^{2m-4}\left(\frac q{Az},\frac q{Bz},\frac q{Cz},\frac q{Dz}\right)_{m-1}
$$
$$
=\sum_{0\le i,j,k,l\le m-1} q^{T_{m-1-i}+T_{m-1-j}+T_{m-1-k}+T_{m-1-l}+m}
$$
\begin{equation} 
\cdot A^iB^jC^kD^l(-z)^{i+j+k+l-2m}\qBin{m-1}{i} \qBin{m-1}{j} \qBin{m-1}{k} \qBin{m-1}{l}. 
\label{eq:6.7}
\end{equation}
Substituting the above expression in the right hand side of \eqn{6.5}, we obtain
$$
S(m)-(1+Aq^m)(1+Bq^m)(1+Cq^m)(1+Dq^m)S(m-1)=
$$
$$
[z^0]\left\{\sum_{0\le i,j,k,l\le m-1}R(A,B,C,D,z) q^{T_{m-1-i}+T_{m-1-j}+T_{m-1-k}+T_{m-1-l}+m}\right.
$$
\begin{equation} 
\cdot  
A^iB^jC^kD^l(z+q^{2m})(-z)^{i+j+k+l-2m}
\qBin{m-1}{i} \qBin{m-1}{j} \qBin{m-1}{k} \qBin{m-1}{l}\bigg\}.
\label{eq:6.8}
\end{equation}
If now the symmetric expansion for $R(A,B,C,D,z)$ in \eqn{6.4} is substituted 
in \eqn{6.8}, then after rearrangement this yields the following recurrence for $S(m)$:
$$
S(m)-(1+Aq^m)(1+Bq^m)(1+Cq^m)(1+Dq^m)S(m-1)=
$$
$$
q^{2m}\mathcal{L}_2(m)-\mathcal{L}_1(m)-(A+B+C+D) \sum_{L\ge 1}F^L
\left\{q^{2m}\mathcal{L}_{1-2L}(m)-\mathcal{L}_{-2L}(m)\right\}
$$
$$
+(ABC+ABD+ACD+BCD)\sum_{L\ge0}F^L
\left\{q^{2m}\mathcal{L}_{-2L}(m)-\mathcal{L}_{-2L-1}(m)\right\}
$$
\begin{equation} 
+F \left\{q^{2m}\mathcal{L}_{-1}(m)-\mathcal{L}_{-2}(m)\right\},
\label{eq:6.9}
\end{equation}
where
$$
\mathcal{L}_n(m)=\sum_{i+j+k+l=2(m-1)+n} q^{T_{m-1-i}+T_{m-1-j}+ 
T_{m-1-k}+T_{m-1-l}+m}A^iB^jC^kD^l 
$$
\begin{equation}
\cdot\qBin{m-1}{i} \qBin{m-1}{j} \qBin{m-1}{k} \qBin{m-1}{l}.
\label{eq:6.10}
\end{equation}
The triangular numbers satisfy the relations
\begin{equation}
T_n=T_{-n-1}\text{ and } T_{m+n}=T_m+T_n+mn
\label{eq:6.11}
\end{equation}
for all integers $m$ and $n$. Using \eqn{6.11} we can rewrite \eqn{6.10}  as
$$
\mathcal{L}_n(m)=q^{-m(n-1)}\sum_{i+j+k+l=2(m-1)+n} q^{T_i+T_j+T_k+T_l}A^iB^jC^kD^l
$$
\begin{equation}
\cdot\qBin{m-1}{i} \qBin{m-1}{j} \qBin{m-1}{k} \qBin{m-1}{l}.
\label{eq:6.12}
\end{equation}

Next we will derive a recurrence relation for the right hand side of 
\eqn{6.1}. To this end we define
\begin{equation}
\sigma_0(m)=\sum_{2m+1\le i+j+k+l} q^{T_i+T_j+T_k+T_l}A^iB^jC^kD^l 
\qBin{m}{i} \qBin{m}{j} \qBin{m}{k} \qBin{m}{l}.
\label{eq:6.13}
\end{equation}
The $q-$binomial coefficients satisfy two recurrence relations, namely,
\begin{equation}
\qBin{m}{i}= \qBin{m-1}{i} +q^{m-i} \qBin{m-1}{i-1}
\label{eq:6.14}
\end{equation}
and
\begin{equation}
\qBin{m}{i}= q^i \qBin{m-1}{i} + \qBin{m-1}{i-1}.
\label{eq:6.15}
\end{equation}
If we use the first recurrence \eqn{6.14} to expand all $q$-binomials in 
\eqn{6.13} we get
\begin{align}
&\sigma_0(m)-(1+Aq^m)(1+Bq^m)(1+Cq^m)(1+Dq^m)\sigma_0(m-1)\nonumber\\
&=-\mathcal{L}_1(m)-q^m\mathcal{L}_2(m)-q^m(A+B+C+D)\mathcal{L}_1(m)\label{eq:6.16}\\
&+q^{2m}(ABC+ABD+ACD+BCD)\mathcal{L}_0(m)+q^{3m}F\mathcal{L}_0(m)+q^{2m} 
F \mathcal{L}_{-1}(m),\nonumber
\end{align}
where $\mathcal{L}_n(m)$ is as in \eqn{6.12}.

Now define
\begin{equation}
\sigma_1(m)=\sum_{L\ge1}q^{2L(m+1)}F^L\sum_{i+j+k+l=2m-2L} 
q^{T_i+T_j+T_k+T_l}A^iB^jC^kD^l 
\qBin{m}{i} \qBin{m}{j} \qBin{m}{k} \qBin{m}{l}.
\label{eq:6.17}
\end{equation}
Next, we expand all the $q-$binomial coefficients in \eqn{6.17} using the second 
recurrence \eqn{6.15}. This yields
\begin{align}
&\sigma_1(m)-(1+Aq^m)(1+Bq^m)(1+Cq^m)(1+Dq^m)\sigma_1(m-1)\nonumber\\
&=q^{3m}F\mathcal{L}_0(m)+(A+B+C+D) \sum_{L\ge1}F^L 
\{q^{2m}\mathcal{L}_{1-2L}(m)-\mathcal{L}_{-2L}(m)\}\nonumber\\
&-(ABC+ABD+ACD+BCD)\sum_{L\ge1}F^L\{q^{2m}\mathcal{L}_{-2L}(m)-\mathcal{L}_{-2L-1}(m)\}\nonumber\\
&-q^{-m}F\mathcal{L}_{-2}(m).
\label{eq:6.18}
\end{align}
The difference
$$
\sigma(m)=\sigma_0(m)-\sigma_1(m)
$$
is precisely the right hand side of \eqn{6.1}. Combining \eqn{6.9}, \eqn{6.16} and \eqn{6.18} 
we find that
\begin{align}
&(S(m)-\sigma(m))-(1+Aq^m)(1+Bq^m)(1+Cq^m)(1+Dq^m)(S(m-1)-\sigma(m-1))\nonumber\\
&=q^m(A+B+C+D)\mathcal{L}_1(m)-(ABC+ABD+ACD+BCD)\mathcal{L}_{-1}(m)\nonumber\\
&+q^m(1+q^m)\mathcal{L}_2(m)-F q^{-m}(1+q^m)\mathcal{L}_{-2}(m).
\label{eq:6.19}
\end{align}
If we show that the right hand side of \eqn{6.19} is $0$, then it follows that $S(m)$ 
and $\sigma(m)$ satisfy the same recurrence.  So this is what we will do next.

To demonstrate that the right hand side of \eqn{6.19} is identically $0$, we 
need to establish that for each $i,j,k,l$, the coefficient of 
$A^iB^jC^kD^l$ is $0$.  Using  \eqn{6.12} and collecting the coefficient of 
$A^iB^jC^kD^l$ in \eqn{6.19} we see that it is sufficient to prove that for
$i+j+k+l=2m$
\begin{align}
&q^m\bigg(q^{T_{i-1}+T_j+T_k+T_l}
\qBin{m-1}{i-1} \qBin{m-1}{j} \qBin{m-1}{k} \qBin{m-1}{l}+\dots\bigg)\nonumber\\
&-q^{2m}\bigg(q^{T_{i-1}+T_{j-1}+T_{k-1}+T_l}
\qBin{m-1}{i-1} \qBin{m-1}{j-1} \qBin{m-1}{k-1} \qBin{m-1}{l}+\dots\bigg)\nonumber\\
+&(1+q^m)\mbox { }q^{T_i+T_j+T_k+T_l}
\qBin{m-1}{i} \qBin{m-1}{j} \qBin{m-1}{k} \qBin{m-1}{l}\nonumber\\
-&q^{2m}(1+q^m)\mbox { }q^{T_{i-1}+T_{j-1}+T_{k-1}+T_{l-1}}
\qBin{m-1}{i-1} \qBin{m-1}{j-1} \qBin{m-1}{k-1} \qBin{m-1}{l-1}\nonumber\\
&=0,
\label{eq:6.20}
\end{align}
where the dots inside the first parenthesis on the left indicate that there 
are three other similar terms in which $i$ is successively interchanged with
$j,k,l$, and the dots in the second parenthesis on the left indicate that 
there are three other terms in which $l$ is succesively interchanged with 
$k,j$ and $i$.

To help us prove \eqn{6.20} we transform the left side using
\begin{equation}
\qBin{m-1}{i-1}=\begin{cases} \frac{1-q^i}{1-q^{m-i}} 
\qBin{m-1}{i}, & \text{ if } m\ne i, \\ 
1, & \text{ if } m=i.
\end{cases}
\label{eq:6.21}
\end{equation}
Since the right hand side of \eqn{6.21} is defined in a piecewise fashion, one 
 needs to consider the following three cases.

\noindent
\underbar{Case A}: $i+j+k+l=2m$ and exactly two of $i,j,k,l$ equal $m$.

Owing to the symmetry in $i,j,k,l$, we may assume $i=j=m$.  Thus $k=l=0$. 
Obviously, each of the $q$-binomial products in \eqn{6.20} is $0$.  
Thus we get no contribution in this case.

\noindent
\underbar{Case B}: $i+j+k+l=2m$ and exactly one of $i,j,k,l$ equals $m$.

Again by symmetry in $i,j,k,l$, we assume $i=m$. So $j+k+l=m$.
In this case the penultimate term in \eqn{6.20} is clearly $0$.  Also in the expression in the 
first parenthesis, 
all three products after the first are $0$ because of the factor 
$\qBin{m-1}{i}$ with $i=m$. Thus using \eqn{6.21} we may rewrite the left hand side of 
\eqn{6.20} in the form
\begin{equation}
q^{T_m+T_j+T_k+T_l}\rho_1(j,k,l,q) \qBin{m-1}{j} \qBin{m-1}{k} \qBin{m-1}{l},
\label{eq:6.22}
\end{equation}
where $j+k+l=m$, and
$$
\rho_1(j,k,l,q)=1-\frac{q^{m-j-k}(1-q^j)(1-q^k)}{(1-q^{m-j})(1-q^{m-k})}-
\frac{q^{m-j-l}(1-q^j)(1-q^l)}{(1-q^{m-j})(1-q^{m-l})}
$$
\begin{equation}
-\frac{q^{m-k-l}(1-q^k)(1-q^l)}{(1-q^{m-k})(1-q^{m-l})}
-\frac{q^{m-j-k-l}(1+q^m)(1-q^j)(1-q^k)(1-q^l)}{(1-q^{m-j})(1-q^{m-k})(1-q^{m-l})}.
\label{eq:6.23}
\end{equation}
It turns out quite miraculously that upon simplification
\begin{equation}
\rho_1(j,k,l,q)=(q^{j+k+l}-q^m)\frac{(1-q^m)^2}{(q^j-q^m)(q^k-q^m)(q^l-q^m)},
\label{eq:6.24}
\end{equation}
and therefore 
\begin{equation}
\rho_1(j,k,l,q)=0 \text{, if }j+k+l=m. 
\label{eq:6.25}
\end{equation}
Thus \eqn{6.25} implies that the expression in \eqn{6.22} has
value $0$, and so there is no contribution from Case B either.

\noindent
\underbar{Case C}: $i+j+k+l=2m$ and none of $i,j,k,l$, equals $m$.

In this case, again using \eqn{6.21}, we may write the expression on the left in \eqn{6.20} as
\begin{equation}
q^{T_i+T_j+T_k+T_l}\rho_2(i,j,k,l,q) 
\qBin{m-1}{i} \qBin{m-1}{j} \qBin{m-1}{k} \qBin{m-1}{l},
\label{eq:6.26}
\end{equation}
where $i+j+k+l=2m$, and
$$
\rho_2(i,j,k,l,q)=\left\{\frac{q^{m-i}(1-q^i)}{1-q^{m-i}}+\dots\right\}-\left\{
\frac{q^{2m-i-j-k}(1-q^i)(1-q^j)(1-q^k)}{(1-q^{m-i}) 
(1-q^{m-j})(1-q^{m-k})}+\dots\right\}
$$
\begin{equation}
+(1+q^m)\left\{1-\frac{(1-q^i)(1-q^j)(1-q^k)(1-q^l)}{(1-q^{m-i})
(1-q^{m-j})(1-q^{m-k})(1-q^{m-l})}\right\}.
\label{eq:6.27}
\end{equation}
In \eqn{6.27} the dots inside the first parenthesis indicate that there are three 
other similar summands with $i$ interchanged with $j,k$, and $l$, in 
succession, and the dots in the second parenthesis indicate that there are 
three other similar summands in which the absence of $l$ is replaced 
successively by the absence of $k,j$, and $i$. If the expression in \eqn{6.27} 
is simplified, it turns out quite miraculously (for a second time) that
\begin{equation}
\rho_2(i,j,k,l,q)=
\label{eq:6.28}
\end{equation}
$$
(q^{2m}-q^{i+j+k+l})\frac{\{(1+q^m)(1-q^i)(1-q^j)(1-q^k)(1-q^l)-(1-q^m)^3\}}
{(q^m-q^i)(q^m-q^j)(q^m-q^k)(q^m-q^l)},
$$
and therefore
\begin{equation}
\rho_2(i,j,k,l,q)=0 \text{, if }i+j+k+l=2m. 
\label{eq:6.29}
\end{equation}
Thus there is no contribution from Case C either. 
Hence, \eqn{6.20} is true.
Clearly,  \eqn{6.19} and \eqn{6.20} together imply that
\begin{equation}
S(m)-\sigma(m)=(1+Aq^m)(1+Bq^m)(1+Cq^m)(1+Dq^m)(S(m-1)-\sigma(m-1)).
\label{eq:6.30}
\end{equation}
Using the definitions \eqn{6.1}, \eqn{6.13} and \eqn{6.17} it is easy to verify
that
\begin{equation}
S(0)=\sigma(0)=0.
\label{eq:6.31}
\end{equation}
From \eqn{6.30} and \eqn{6.31} it follows that
\begin{equation}
S(m)=\sigma(m)\text{ for all } m\ge 0.
\label{eq:6.32}
\end{equation}
Finally, observe that \eqn{6.13} and \eqn{6.17} imply
\begin{equation}
\begin{cases}
\left[q^n\right]\sigma_0(m)=0, & \text{ if } n<T_{m+1},\\
\left[q^n\right]\sigma_1(m)=0, & \text{ if } n<T_{m+1}+\lfloor\frac{m+2}{2}\rfloor, 
\end{cases}
\label{eq:6.33}
\end{equation}
and therefore
\begin{equation}
\lim\limits_{m\rightarrow\infty} \sigma_0(m)=\lim\limits_{m\rightarrow\infty} 
\sigma_1(m)=0
\label{eq:6.34}
\end{equation}
Since $\sigma(m)=\sigma_0(m)-\sigma_1(m)$, \eqn{6.34} yields
\begin{equation}
\lim\limits_{m\rightarrow\infty}\sigma(m)=0.
\label{eq:6.35}
\end{equation}
This combined with \eqn{6.32} shows
\begin{equation}
\lim\limits_{m\rightarrow\infty} S(m)=\lim\limits_{m\rightarrow\infty}\sigma(m)=0
\label{eq:6.36}
\end{equation}
and this completes the proof of \eqn{5.18} and, as a result, of the key identity \eqn{1.2}.

\bigskip

\section{\bf Consequences and prospects}

\medskip

In order to obtain G\"ollnitz's theorem from the identity 
\eqn{ext.2}, the transformations \eqn{2.16} were used. The dilation $q\mapsto q^6$ in \eqn{2.16}
is the minimal one required to convert the six colors (three primary and three secondary) in 
Theorem~\ref{thm:5} to non--overlapping residue classes. If a dilation $q\mapsto q^M$ with
$M<6$ is used, then the six colors will lead to some overlap in the residue classes, that is 
to say, certain residue classes mod $M$ will represent integers in more than one of the six 
colors.  This will involve a discussion of partitions with weights attached, these weights 
being polynomials in $A,B$, and $C$. A general theory of weighted partition identities and 
its uses can be found in Alladi \cite{2}. In addition Alladi \cite{1}, \cite{5}, and more 
recently Alladi and Berkovich \cite{13}, \cite{15} discuss weighted partition reformulations
of G\"ollnitz's theorem obtained from dilations $q\mapsto q^M$ with $M<6$, and some of their important
consequences. In some instances, values of the parameters $A,B,C$ can be chosen so that there 
is a major collapse of the weights.  This has led to derivations of various equivalent forms of
Jacobi's celebrated triple product identity \eqn{5.8}. 

To get Theorem~\ref{thm:3} from 
Theorem~\ref{thm:6} we used the dilation $q\mapsto q^{15}$ in \eqn{1.3}. This dilation led to
$11$ distinct residue classes mod $15$ being used in Theorem 1.  Since Theorem~\ref{thm:6} 
involves $11$ colors, any dilation $q\mapsto q^M$ with $M<11$ will lead to an overlap in 
residue classes and therefore to weighted partition identities.  It would be worthwhile to
discuss weighted partition identities that emerge from Theorem~\ref{thm:6} in a manner similar
to what was done in \cite{1}, \cite{5}, \cite{13}, \cite{15}.  Another very interesting question
is whether there exists a  quadruple product extension of \eqn{5.8}.
More precisely, observing that \eqn{5.8} contains a lacunary expansion of a triple product
with one free parameter $z$, we ask whether there is an expansion of a quadruple product 
with two free parameters which reduces to \eqn{5.8} when one of the parameters is set equal 
to $0$.  Recently, by setting $D=-1$ in \eqn{1.2}, Alladi and Berkovich \cite{14}
succeeded in reformulating Theorem~\ref{thm:6} as a weighted partition identity. It 
turned out that the choice
$$
C=B^{-1}
$$
(motivated by $z\cdot z^{-1}=1$ in \eqn{5.8}) led to a significant collapse in the 
weights, and to the following quadruple product extension of \eqn{5.8}:\\
{\bf Proposition:}
\begin{align}
&\sum_{k,l,m,n\ge0} q^{T_l}\left(\frac{B^{l+1}+B^{-l}}{B+1}\right) 
\frac{q^{4T_n+2T_m+k^2+2n(l+m)+m(l-k)}(-1)^nA^{n+m}B^{2k-m}}{(q)_n(q)_k(q)_{m-k}} \nonumber\\
+&\sum_{\substack{k,m,n\ge0\\ l\ge1}}A^lq^{T_l} 
\frac{q^{4T_n+2T_m+k^2+2n(l+m)+m(l-k)-n}(-1)^nA^{n+m}B^{2k-m}}{(q)_n(q)_k(q)_{m-k}} \nonumber
\end{align}
\begin{equation}
=(-Aq,-Bq,-B^{-1}q,q)_\infty.
\label{eq:7.2}
\end{equation}

If we put $A=0$ in \eqn{7.2}, it forces the collapse of the entire second sum 
on the left; also the contribution from the first sum in \eqn{7.2} is only when 
$k=m=n=0$. So with $A=0$, \eqn{7.2} reduces to
\begin{equation}
\sum^\infty_{l=0} 
q^{T_l}\left(\frac{B^{l+1}+B^{-l}}{B+1}\right)=(-Bq,-B^{-1}q,q)_\infty,
\label{eq:7.3}
\end{equation}
which is the same as \eqn{5.8} with $B=z$. Being a theta function identity, 
\eqn{7.3} is very fundamental. It would be interesting to investigate 
applications of the more general identity \eqn{7.2}. The derivation of \eqn{7.2} as 
a special case of the key identity \eqn{1.2} involves some nontrivial combinatorial analysis 
and will be presented elsewhere \cite{14}.

In the course of studying vertex operators in the theory of affine Lie algebras, 
Capparelli \cite{23} was led to the following theorem: {\it The number of partitions 
$C^*(n)$ of $n$ into parts $\equiv\pm2$, $\pm3$ (mod $12$) equals the number 
of partitions $D(n)$ of $n$ into parts $>1$, and differing by $\ge2$, where 
the difference is $\ge4$ unless consecutive parts are multiples of $3$ or add up to 
a multiple of $6$}.

The generating function of $C^*(n)$ is
\begin{equation}
\sum^\infty_{n=0} C^*(n)q^n=\frac 1{(q^2,q^3,q^9,q^{10};q^{12})_\infty}.
\label{eq:7.4}
\end{equation}
Alladi, Andrews, and Gordon \cite{9} showed that by rewriting the product in \eqn{7.4} as
$$
(-q^2,-q^3,-q^4,-q^6;q^6)_\infty,
$$
it is possible to get a two parameter refinement of Capparelli's 
theorem.  More precisely, if $C(n;i,j)$ is defined by
\begin{equation}
\sum_{n,i,j}C(n;i,j)A^iB^jq^n=(-Aq^2,-q^3,-Bq^4,-q^6;q^6)_\infty,
\label{eq:7.5}
\end{equation}
then there is a two parameter refinement of Capparelli's theorem in 
which  $C(n,i,j)$ is equated to $D(n;i,j)$, defined as $D(n)$ above subject to the additional 
restrictions that there are precisely $i$ parts $\equiv 2\pmod{3}$ and $j$ parts $\equiv 1\pmod{3}$.  In \cite{9} a combinatorial bijective proof of this two parameter 
refined theorem is given.

Notice that the transformations
\begin{equation}
\begin{cases}
\text{(dilation) } q\mapsto q^6, \\
\text{(translations) }A\mapsto Aq^{-4},B\mapsto Bq^{-2},C\mapsto 
Cq^{-3},D\mapsto D,
\end{cases}
\label{eq:7.6}
\end{equation}
in \eqn{1.2} lead to the product
\begin{equation}
(-Aq^2,-Cq^3,-Bq^4,-Dq^6;q^6)_\infty,
\label{eq:7.7}
\end{equation}
which is a four parameter refinement of the product in \eqn{7.5}. Since 
$q\mapsto q^6$ is a dilation smaller than $q\mapsto q^{11}$, in view of 
certain remarks made above, \eqn{1.2} and Theorem~\ref{thm:6} will lead to a weighted 
partition theorem in which partitions with parts differing by $\ge6$ and 
with weights (polynomials in $A,B,C,D$) attached are equated to the 
partition function whose generating function is the product in \eqn{7.7}. The 
gap conditions that would emerge from this construction are different from those defining
 $D(n)$ in Capparelli's theorem. It would be worthwhile to investigate whether the 
gap conditions for our weighted partitions are connected to those defining 
$D(n)$, and whether our companion of Capparelli's theorem has implications in the theory of
affine Lie algebras.

Motivated by studies in the representation theory of symmetric groups, 
Andrews, Bessenrodt, and Olsson \cite{21} considered partitions whose generating 
function is the product
\begin{equation}
(-\frac ax q,-\frac ay q^2,-byq^3,-bxq^4;q^5)_\infty
\label{eq:7.8}
\end{equation}
and related these to partitions satisfying certain difference conditions. 
We may view the product in \eqn{7.8} as one of the form
\begin{equation}
(-Aq,-Bq^2,-Cq^3,-Dq^4;q^5)_\infty,
\label{eq:7.9}
\end{equation}
where there is the relation
$$
AD=BC
$$
among the parameters. Thus there are only three free parameters in \eqn{7.8}. 
Clearly the transformations
\begin{equation}
\begin{cases}
\text{(dilation) } q\mapsto q^5, \\
\text{(translations) } A\mapsto Aq^{-4}, B\mapsto Bq^{-3}, C\mapsto 
Cq^{-2}, D\mapsto Dq^{-1}
\end{cases}
\label{eq:7.10}
\end{equation}
convert the product in \eqn{1.2} to that in \eqn{7.9}. Thus in this case, Theorem~\ref{thm:6} 
will yield a result in which a four parameter weighted count of certain 
partitions into parts differing by $\ge5$ are connected to partitions whose 
generating function is the product in \eqn{7.9}. As in the case of Capparelli's theorem, the 
gap conditions we would get are different from those considered in \cite{21}, but perhaps there
are connections  worth exploring.

In the last decade, many new generalizations of the Rogers--Ramanujan 
identities were discovered and proved by McCoy and collaborators (see \cite{22} 
for a review and references) using the so called themodynamic Bethe ansatz 
(TBA) techniques.  It would be highly desirable to find a TBA 
interpretation of the new identity \eqn{1.7}. Such an intepretation, besides 
being of substantial interest in physics, may provide insight into how to 
extend Theorem~\ref{thm:6} to the one with five or more primary colors.

Recently Alladi and Berkovich \cite{12} obtained a bounded version of \eqn{2.15}:
\begin{align}
&\sum_{\substack{i=a+ab+ac\\ j=b+ab+bc\\ k=c+ac+bc}} 
q^{T_s+T_{ab}+T_{ac}+T_{bc-1}} \nonumber\\
&\cdot\Bigg\{ q^{bc} \qBin{L-s+a}{a} \qBin{L-s+b}{b}
   \qBin{L-s+c}{c} \qBin{L-s}{ab} \qBin{L-s}{ac} \qBin{L-s}{bc} \notag\\
   &+ \qBin{L-s+a-1}{a-1} \qBin{L-s+b}{b}
   \qBin{L-s+c}{c} \qBin{L-s}{ab} \qBin{L-s}{ac} \qBin{L-s}{bc-1} \Bigg\} \notag\\
   &= q^{T_i+T_j+T_k}
   \qBin{L-i}{j} \qBin{L-j}{k} \qBin{L-k}{i},
\label{eq:7.11}
\end{align}
where $s=a+b+c+ab+ac+bc$. 
Actually in \cite{12}, a doubly bounded version of \eqn{2.15} was established. Subsequently, 
the full triply bounded refinement of \eqn{2.15} was found and proven by Berkovich and Riese 
in \cite{BR}. 

Alladi and Berkovich \cite{12} noticed that \eqn{7.11} can 
be interpreted combinatorially to yield a refinement of Theorem~\ref{thm:5} with 
bounds on the sizes of the parts. By considering dilations 
$q\mapsto q^M$ with $M<6$, it can be shown that the bounded G\"ollnitz 
theorem leads to new bounded versions of fundamental $q-$series identities including
\eqn{5.8}.  A detailed treatment of some of these important applications may be found 
in \cite{15}.
In view of \eqn{7.11} we may ask whether there exists a bounded version of the 
new four parameter identity \eqn{1.7}. 

Finally, it is important to investigate what 
kind of $q-$hypergeometric transformation formulae hidden behind \eqn{1.7} and \eqn{ins2}.
This question as well as those 
mentioned above indicate that \eqn{1.7} opens up several exciting avenues for further 
exploration.

\subsection*{Acknowledgement}

We are grateful to Frank Garvan, Barry McCoy, Peter Paule, Axel Riese and Doron Zeilberger
for their kind interest and comments. We thank the referee for a careful reading of the manuscript 
and helpful suggestions.

\end{document}